\documentclass[reqno]{amsart}
\usepackage[english]{babel}
\usepackage{amscd,amssymb,amsmath,amsfonts,latexsym,amsthm}
\usepackage{inputenc}
\usepackage{graphicx}
\usepackage{color}

\newtheorem{thm}{Theorem}[section]

\newtheorem{prop}[thm]{Proposition}
\newtheorem{defn}[thm]{Definition}

\newtheorem{rem}[thm]{Remark}

\numberwithin{equation}{section}
\DeclareFontFamily{OT1}{pzc}{}
\DeclareFontShape{OT1}{pzc}{m}{it}{<-> s * [1.10] pzcmi7t}{}
\DeclareMathAlphabet{\mathpzc}{OT1}{pzc}{m}{it}

\title{Biderivations of some classes of solvable Leibniz algebras}

\author[B.B.Yusupov]{Bakhtiyor Yusupov$^{1,2}$}

\author[D.E.Jumaniyozov]{Doston Jumaniyozov$^{1,3}$}

\author[M.E.Azizov]{Majidkhon Azizov$^1$}

\email{baxtiyor\_yusupov\_93@mail.ru, dostondzhuma@gmail.com, azizovmajidkhan@gmail.com}

\address{$^1$ V.I.Romanovskiy Institute of Mathematics Uzbekistan Academy of Sciences. 9,  University street, 100174, Tashkent, Uzbekistan}

\address{$^2$ Department of Algebra and Mathematical engineering, Urgench State University, H. Alimdjan street, 14, Urgench 220100, Uzbekistan}

\address{$^3$ National University of Uzbekistan, University street, 4, Olmazor district, Tashkent, 100174, Uzbekistan}

\begin{document}
\maketitle
\noindent {\bf Abstract.}
{\it In this work, we investigate anti-derivations and biderivation of Leibniz algebras. We describe general form of anti-derivations and biderivations on null-filiform and filiform Leibniz algebras. Moreover, we show how to construct Leibniz algebras, while using these biderivations. We describe general form of anti-derivations and biderivations on solvable Leibniz algebras with  null-filiform and filiform nilradicals. Interesting fact that, any biderivations of solvable Leibniz algebras with  null-filiform and filiform nilradicals are inner biderivations.
}

\bigskip



\section{Introduction}

Leibniz algebras, first introduced by Blokh in 1965 as a broader framework than Lie algebras \cite{blox}, were later thoroughly investigated by Loday in 1993 \cite{loday1993version}. The problem of classification of finite-dimensional Leibniz algebras is fundamental and a very complicated problem. During the last three decades, significant attention has been given to the study of Leibniz algebras, with many publications dedicated to their investigation \cite{Ayupov1, Ayupov2, Ayupov5, Ayubook}. The analogue of the Levi-Malcev decomposition for Leibniz algebras was proved by D.W.Barnes \cite{Bar}, that asserts that any Leibniz algebra decomposes into a semidirect sum of its solvable radical and a semisimple Lie algebra. The semisimple part of a Leibniz algebra can be described in terms of its simple Lie ideals. Consequently, the main challenge in classifying finite-dimensional Leibniz algebras reduces to the study of solvable Leibniz algebras. In this direction, significant progress has been made in the classification of solvable Leibniz algebras with prescribed nilradicals \cite{CLOK12, CLOK13, ladra, KhS}.

Among the fundamental concepts in the study of Leibniz algebras are derivations, anti-derivations, automorphism and local derivations, which describe how the structure of the algebra evolves under certain linear maps \cite{AdaYus, Ayupov5, Bar, YusBS}. A natural progression in this line of inquiry involves the examination of biderivations, defined as ordered pairs comprising a derivation and an antiderivation that fulfill certain compatibility relations. The concept of biderivations was originally formulated by Brešar in \cite{bres} within the context of associative algebras and has since been generalized to encompass broader algebraic frameworks, including both Lie and Leibniz algebras.
Especially, Bre\u{s}ar et al. introduced the notion of biderivation of rings in \cite{B4}, and they showed that all biderivations of noncommutative prime rings are inner. In \cite{W2}, the authors proved that all skew-symmetric biderivations of finite dimensional simple Lie algebras over an algebraically closed field of characteristic zero are inner.
Originally, the notion of biderivation for a Leibniz algebra $L$ was first introduced by Loday in 1993 as a pair of linear maps $(d,D),$ where $d$ is a derivation and $D$ is an anti-derivation satisfying the additional relation $[x,d(y)]=[x,D(y)]$ for all $x,y\in L.$
In \cite{M} given a complete classification of the Leibniz algebras of biderivations of right Leibniz algebras of dimension up to three over a field $F$, with $char(F) \neq 2$ and described the main properties of such class of Leibniz algebras. In \cite{ZM} focused on the biderivations of 4-dimensional nilpotent complex Leibniz algebras. Using the existing classification of these algebras, Zahari was developed algorithms to compute derivations, antiderivations, and biderivations as pairs of matrices with respect to a fixed basis.

The investigation of biderivations offers significant insight into the underlying structural characteristics of algebras. In the context of Leibniz algebras, biderivations serve as tools for identifying symmetries, invariant subspaces, and other essential algebraic features. Although the theory of derivations and inner derivations has been extensively developed for Lie and associative algebras, the computational study of biderivations within Leibniz algebras remains relatively limited, particularly in the context of 4-dimensional algebras, despite their presence in existing classification results.

This study focuses on anti-derivations and biderivations in solvable Leibniz algebras, with particular attention to the concept of biderivations as originally introduced by Loday in \cite{loday1993version}. The structure of the paper is as follows. Section 2 presents the fundamental concepts and definitions related to Leibniz algebras. In Section 3, we characterize the general forms of anti-derivations and biderivations for null-filiform and filiform Leibniz algebras. Additionally, we construct the associated Leibniz algebras of biderivations for these cases. Section 4 is devoted to the study of anti-derivations and biderivations in solvable Leibniz algebras whose nilradicals are either null-filiform or filiform. A noteworthy result established in this section is that all biderivations in such solvable Leibniz algebras turn out to be inner.

\section{Preliminaries}

In this section we recall the definitions of \emph{derivation}, \emph{anti-derivation} and \emph{biderivation} for right Leibniz algebras and we show an example of computation of the biderivations algebra.

\begin{defn}
	Let $L$ be a Leibniz algebra over $\mathbb{F}$. A \emph{derivation} of $L$ is a linear map $d\colon L\rightarrow L$ such that
	$$
	d(\left[x,y\right])=\left[d(x),y\right]+\left[x,d(y)\right],\;\; \forall x,y\in L.
	$$
\end{defn}

The right multiplications are particular derivations called \emph{inner} \emph{derivations} and an equivalent way to define a right Leibniz algebra is to saying that the (right) adjoint map $\operatorname{ad}_x=\left[-,x\right]$ is a derivation, for every $x\in L$. Meanwhile, for a right Leibniz algebra, the left adjoint maps are not derivations in general.

With the usual bracket \hbox{$\left[d_1,d_2\right]=d_1\circ d_2 - d_2\circ d_1$}, the set $\operatorname{Der}(L)$ is a Lie algebra and the set  $\operatorname{Inn}(L)$ of all inner derivations of $L$ is an ideal of $\operatorname{Der}(L)$. Furthermore, $\operatorname{Aut}(L)$ is a Lie group and the associated Lie algebra is $\operatorname{Der}(L)$.\\

The definitions of \emph{anti-derivation} and \emph{biderivation} for a Leibniz algebra were first given by J.-L.\ Loday in \cite{loday1993version}.

\begin{defn}
	An \emph{anti-derivation} of a right Leibniz algebra $L$ is a linear map $D\colon L \rightarrow L$ such that
	$$
	D([x,y])=[D(x),y]-[D(y),x], \; \; \forall x,y \in L.
	$$
\end{defn}
\noindent
For a left Leibniz algebra we have to ask that
$$
D([x,y])=[x,D(y)]-[y,D(x)], \; \; \forall x,y \in L.
$$
We observe that in the case of Lie algebras, there is no difference between a derivation and an anti-derivation. Moreover one can check that, for every $x \in L$, the left adjoint map
\begin{equation}\label{Antider}
    \operatorname{Ad}_x=[x,-]
\end{equation}
defines an anti-derivation.

\begin{rem}
	The set of anti-derivations of a Leibniz algebra $L$ has a $\operatorname{Der}(L)-$module structure with the action
	$$
	d \cdot D := [D,d]= D \circ d - d \circ D,
	$$
	for every $d \in \operatorname{Der}(L)$ and for every anti-derivation $D$.
\end{rem}

\begin{rem}
	Let $L$ be a Leibniz algebra.  $D\colon L \rightarrow L$ be an anti-derivation. Then, for every $x \in L$, we have
	$$
	D([x,x])=[D(x),x]-[D(x),x]=0,
	$$
	thus $D(L)=0$.
\end{rem}

\begin{defn}
	Let $L$ be a right Leibniz algebra. A \emph{biderivation} of $L$ is a pair
	$$
	(d,D)
	$$
	where $d$ is a derivation and $D$ is an anti-derivation, such that
	\begin{equation}\label{1}
		[x,d(y)]=[x,D(y)], \; \; \forall x,y \in L.
	\end{equation}
\end{defn}

The set of all biderivations of $L$, denoted by $\operatorname{Bider}(L)$, has a Leibniz algebra structure with the bracket
\begin{align}\label{bider}
[(d,D),(d',D')]=(d \circ d' - d' \circ d, D \circ d' - d' \circ D), \; \; \forall (d,D),(d,D') \in \operatorname{Bider}(L),
\end{align}
and it is possible to define a Leibniz algebra morphism
$$
L \rightarrow \operatorname{Bider}(L)
$$
by
$$
x \mapsto (-\operatorname{ad}_x, \operatorname{Ad}_x), \; \; \forall x \in L.
$$
The pair $(-\operatorname{ad}_x, \operatorname{Ad}_x)$ is called \emph{inner biderivation} of $L$ and the set of all inner biderivations forms a Leibniz subalgebra of $\operatorname{Bider}(L)$.

\begin{rem}
It should be noted that in \cite{CKL25}, another notion of biderivation is defined. Namely, a bilinear map $\delta:L\times L\to L$ is called a {\em  biderivation} if for any $x,y,z\in L$, the following identities hold true:
\begin{align*}
\delta([x,y],z)&=[x,\delta(y,z)]-[y,\delta(x,z)],\\
\delta(x, [y,z])&=[\delta(x,y),z]+[y, \delta(x, z)].
\end{align*}
If, furthermore, $\delta(x,y) = \delta(y,x)$ for all $x,y\in L$, we call $\delta$ a symmetric biderivation. Here, we note that if $\delta$ is a biderivation, then the identities above mean that for any $z\in L$ the linear maps $d: L\to L, \ d(x):=\delta(z,x)$ and $D: L\to L, \ D(x):=\delta(x,z)$ are derivation and antiderivation, respectively.
\end{rem}

Let $\mathcal{L}$ be a Leibniz algebra. For a Leibniz algebra $\mathcal{L}$ consider the following central lower and
derived sequences:
$$
\mathcal{L}^1=\mathcal{L},\quad \mathcal{L}^{k+1}=[\mathcal{L}^k,\mathcal{L}^1], \quad k \geq 1,
$$
$$\mathcal{L}^{[1]} =\mathcal{L}, \quad \mathcal{L}^{[s+1]} = [\mathcal{L}^{[s]}, \mathcal{L}^{[s]}], \quad s \geq 1.$$

\begin{defn} A Leibniz algebra $\mathcal{L}$ is called
nilpotent (respectively, solvable), if there exists  $p\in\mathbb N$ $(q\in
\mathbb N)$ such that $\mathcal{L}^p=0$ (respectively, $\mathcal{L}^{[q]}=0$).The minimal number $p$ (respectively, $q$) with such
property is said to be the index of nilpotency (respectively, of solvability) of the algebra $\mathcal{L}$.
\end{defn}

Note that any Leibniz algebra $\mathcal{L}$ contains a unique maximal solvable (resp. nilpotent) ideal, called the radical (resp. nilradical) of the algebra.

\section{Biderivations on nilpotent Leibniz algebras}

It is known that the maximal nilindex for the $n$-dimensional Leibniz algebras is equal to $n+1$ and such Leibniz algebras are called null-filiform. In \cite{Ayupov5} it is proved that up to isomorphism there exists a unique null-filiform Leibniz algebra with basis $\{e_1,e_2,\dots,e_n\}$ and the multiplication table defined as follows:
$$ NF_{n}:[e_{i},e_{1}]=e_{i+1} \ \ \ \ \ 1\leq i \leq n-1.$$

\begin{defn}
An $n$-dimensional Leibniz algebra $\mathcal{L}$ is said to be filiform if $
\dim \mathcal{L}^i=n-i$ for $2\leq i \leq n$.
\end{defn}

In \cite{Ayupov5} it is shown that up to isomorphism there exist two non-Lie naturally graded filiform Leibniz algebras

\begin{itemize}
\item $F_{n}^{1}:  [e_{1},e_{1}]=e_{3}, \ \ \
[e_{i},e_{1}]=e_{i+1}$,\ \ \ \   ${2}\leq{i}\leq{n-1}$
\item $F_{n}^{2}:  [e_{1},e_{1}]=e_{3}, \ \ \ [e_{i},e_{1}]=e_{i+1}$,\ \ \ \   ${3}\leq{i}\leq{n-1}.$
\end{itemize}

\subsection{Biderivations of null-filiform Leibniz algebras}
In this part, we classify the biderivations of null-filiform Leibniz algebras and construct induced Leibniz algebras. Moreover, we explore the relation between the algebras and the Leibniz algebra of its biderivations.
First, we provide descriptions of the derivations and anti-derivations on the given Leibniz algebras. Then we give the classification of the biderivations as well as the Leibniz algebras constructed by them.

\begin{prop}\cite{CLOK12} Any derivation on $NF_n$ has the following form:
    \begin{equation*}
        d(e_i)=i\alpha_1e_i+\sum\limits_{j=i+1}^n\alpha_{j-i+1}e_j.
    \end{equation*}
where $\alpha_i\in\mathbb{C},$ $i\in\{1,2,\dots,n\}.$
\end{prop}

Now, we determine the space of anti-derivations of $NF_n.$

\begin{prop} Any anti-derivation on $NF_n$ has the following form:
    \begin{equation*}
        D(e_1)=\sum\limits_{j=1}^n\beta_je_j,\quad D(e_i)=0,\quad 2\leq i\leq n.
    \end{equation*}
where $\beta_j\in\mathbb{C},$ $j\in\{1,2,\dots,n\}.$
\end{prop}

\begin{proof} Let $D$ be an anti-derivation of $NF_n.$ Set
\[D(e_1)=\sum\limits_{i=1}^n\beta_ie_i.\]
Then by $[e_1,e_1]=e_2$ and applying the identity of anti-derivation, we get the following restriction to $D$:
\begin{equation*}
    \begin{split}
        D(e_2)=D([e_1,e_1])=[D(e_1),e_1]-[D(e_1),e_1]=0.
    \end{split}
\end{equation*}
Applying the same arguments for $[e_{i-1},e_1]=e_i$, we get the following relations:
\begin{equation*}
    \begin{split}
D(e_i)=D([e_{i-1},e_1])=[D(e_{i-1}),e_1]-[D(e_1),e_{i-1}]=0,
    \end{split}
\end{equation*}
where $i \in \lbrace 2, 3, \ldots n \rbrace.$ The remaining products reduce to an identity.
\end{proof}

\begin{thm}
    Any biderivation of $NF_n$ is described as the following form:
\begin{equation}
\left\{ \left(
\begin{array}{cccccccccc}
 \alpha_1    & 0             & \cdots & 0          & 0                  \\
 \alpha_2    & 2\alpha_1   & \cdots & 0          & 0                \\
 \alpha_3    & \alpha_2    & \cdots & 0          & 0                \\
 \cdots        & \cdots        & \cdots & \cdots     & \cdots           \\
 \alpha_{n-1}& \alpha_{n-2}& \cdots & (n-1)\alpha_1 & 0            \\
 \alpha_n    &\alpha_{n-1} & \cdots & \alpha_2 & n\alpha_1        \\
 \end{array}
\right), \
\left(\begin{array}{cccccccccc}
 \alpha_1    & 0             & \cdots & 0          & 0                  \\
 \beta_2    & 0  & \cdots & 0          & 0                \\
 \beta_3    & 0   & \cdots & 0          & 0                \\
 \cdots        & \cdots        & \cdots & \cdots     & \cdots           \\
 \beta_{n-1}& 0& \cdots & 0& 0            \\
 \beta_n    &0 & \cdots & 0 & 0        \\
 \end{array}\right)\right\}
    \end{equation}
\end{thm}
\begin{proof} Let $(d,D)$ be a biderivation of $NF_n,$ where $d\in Der(NF_n)$ and $D$ is an anti-derivation. We first apply the identity \eqref{1} for the elements $x=e_1, \ y=e_1$ and
  substitute the values from the derivation and anti-derivation. Then we have
    $$[e_1,\alpha_1e_1+\alpha_2e_2+\cdots+\alpha_ne_n]=[e_1,\beta_1e_1+\beta_2e_2+\cdots+\beta_ne_n].$$
Comparing the corresponding basis coefficients, we get that $\beta_1=\alpha_1.$ One can verify that the remaining products yields the identity.
    \end{proof}

Now, we describe the Leibniz algebra ${\rm Bider}(NF_n)$.
For this, let us fix the basis $\mathcal{B}$ of ${\rm Bider}(NF_n)$. One can verify that ${\rm Bider}(NF_n)$ is $(2n-1)$-dimensional. Thus, consider $\mathcal{B}=\{X_1,X_2,\dots,X_n,Y_2,Y_3,\dots,Y_{n}\}$, where
\begin{align*}
    X_1&=\left(E_{11}+2E_{22}+\dots+nE_{nn},E_{11}\right)   & &\\
    X_2&=\left(E_{21}+E_{32}+\dots+E_{nn-1},0\right)    &  Y_2=&\left( 0,E_{21} \right)\\
    &\dots & \dots\\
     X_{n-1}&=\left(E_{n-11}+E_{n2},0\right)   & Y_{n-1}=&\left( 0,E_{n-11} \right)\\
     X_{n}&=\left(E_{n1},0\right)   & Y_{n}=&\left( 0,E_{n1} \right),
\end{align*}
where $E_{ij}$ are matrix units.
Then it can be verified by \eqref{bider} that $[X_1,X_1]=[X_i,X_j]=0$ and $[X_1,Y_j]=[X_i,Y_j]=[Y_i,Y_j]=0,$ for $2\leq i,j\leq n$.
For other products we do the following computations:
\begin{align*}
    [X_1,X_k]&=\left([E_{11}+2E_{22}+\dots+nE_{nn},E_{k1}+\dots+E_{nn-k+1}],[E_{11},E_{k1}+\dots+E_{nn-k+1}]\right)\\
    &=\left(E_{k1}+E_{k+12}+\dots+E_{nn-k+1},-E_{k1}\right)=X_k-Y_k, \ \mbox{ for  } \ 2\leq k\leq n,\\
    [X_k,X_1]&=\left([E_{k1}+\dots+E_{nn-k+1},E_{11}+2E_{22}+\dots+nE_{nn}],0\right)\\
    &=-\left(E_{k1}+E_{k+12}+\dots+E_{nn-k+1},0\right)=-X_k, \ \mbox{ for } \ 2\leq k\leq n,\\
    [Y_k,X_1]&=\left(0,[E_{k1},E_{11}+2E_{22}+\dots+nE_{nn}]\right)=-(k-1)\left(0,E_{k 1}\right)=-(k-1)Y_k\\
    [Y_k,X_l]&=\left(0,[E_{k1},E_{l1}+E_{l+12}+\dots+E_{nn-l+1}]\right)=-\left(0,E_{k+l-11}\right)=-Y_{k+l-1} \ \mbox{ for } \ k+l-1\leq n.
\end{align*}

Thus, by summarizing all above, we obtain the following basis products:
\begin{align}\label{table1}
    [X_1,X_k]&=X_k-Y_k, & [X_k,X_1]&=-X_k,\\
    \label{table2}
    [Y_k,X_1]&=-(k-1)Y_k &  [Y_k,X_l]&=-Y_{k+l-1},
\end{align}
where $2\leq k\leq n, \ k+l-1\leq n$ and the omitted products are zero.

\begin{prop}
    ${\rm Bider}(NF_n)$ is a solvable Leibniz algebra of dimension $2n-1$ with the multiplication table given as \eqref{table1} and \eqref{table2}.
\end{prop}

\subsection{Biderivations of filiform Leibniz algebras}
In this subsection, we describe the spaces of biderivations of filiform Leibniz algebras. First of all, we provide a known result regarding the space of derivations of filiform Leibniz algebras $F_n^1$ and $F_n^2.$

\begin{prop}\cite{CLOK13}
 Any derivation of the algebras $F_n^1$ and $F_n^2$ has the following form:

$\bullet$ For algebra $F_n^1:$
 $$
 \begin{pmatrix}
\alpha_1 & 0 & 0 & \cdots & 0\\
\alpha_2 & \alpha_1+\alpha_2 & 0 & \cdots & 0\\
\alpha_3 & \alpha_3 & 2\alpha_1+\alpha_2 & \cdots & 0\\
\vdots & \vdots & \vdots & \ddots & \vdots\\
\alpha_n & \alpha_{n-1} & \alpha_{n-2} & \cdots & (n-1)\alpha_1+\alpha_2
\end{pmatrix},\qquad
 $$

 $\bullet$ For algebra $F_n^2:$
 $$
 \begin{pmatrix}
\alpha_1 & 0 & 0 & \cdots & 0\\
\alpha_2 & \alpha_{n+1} & 0 & \cdots & 0\\
\alpha_3 & 0 & 2\alpha_1 & \cdots & 0\\
\vdots & \vdots & \vdots & \ddots & \vdots\\
\alpha_n & 0 & \alpha_4 & \cdots & (n-1)\alpha_1
\end{pmatrix}.
 $$
where $\alpha_i\in\mathbb{C}.$
\end{prop}

Now, we determine the spaces of anti-derivations of $F_n^1$ and $F_n^2.$
\begin{prop}
Any anti-derivation of the algebras $F_n^1$ and $F_n^2$ has the following form:

    $\bullet$ For algebra $F_n^1:$

    \begin{equation*}
\begin{split}
D_1(e_1)&=\sum\limits_{j=1}^n\beta_je_j,\quad  D_1(e_2)=\beta_{n+1}e_1-\beta_{n+1}e_2+\beta_{n+2}e_n, \quad
D_1(e_i)=0,\quad 3\leq i\leq n.
\end{split}
\end{equation*}

$\bullet$ For algebra $F_n^2:$

\begin{equation*}
\begin{split}
D_2(e_1)&=\sum\limits_{j=1}^n\beta_je_j,\quad  D_2(e_2)=\beta_{n+1}e_2+\beta_{n+2}2 e_n, \quad
D_2(e_i)=0,\quad 3\leq i\leq n.
\end{split}
\end{equation*}
where $\beta_i\in\mathbb{C}.$
\end{prop}
\begin{proof}  We prove the proposition for the algebra $F_n^1$, and for the algebra $F_n^2$ the proof is similar.
    It is easy to see that the basis elements $e_1$ and $e_2$ are the generators of $F_n^1.$ First we define the action of an anti-derivation on these generators, and then applying the identity of anti-derivation we define the action on the other elements. Let $D$ be an anti-derivation of $F_n^1.$ We set
$$D(e_1)=\sum\limits_{i=1}^{n}\beta_{1i}e_i,\quad
D(e_2)=\sum\limits_{i=1}^{n}\beta_{2i}e_i.$$
Now consider the condition of anti-derivation for the elements $e_1$ and $e_2:$
\begin{equation*}
    \begin{split}
            D(e_3)&=D([e_2,e_1])=[D(e_2),e_1]-[D(e_1),e_2]=\sum\limits_{i=1}^{n}\beta_{2i}[e_i,e_1]-\sum\limits_{i=1}^{n}\beta_{1i}[e_i,e_2]=\\
            &=(\beta_{21}+\beta_{22})e_3+\sum\limits_{i=4}^n\beta_{2 \ i-1}e_i.
    \end{split}
\end{equation*}
On the other hand, we have
\begin{equation*}
    \begin{split}
    D(e_3)&=D([e_1,e_1])=[D(e_1),e_1]-[D(e_1),e_1]=0.\\
 \end{split}
\end{equation*}
Comparing coefficients of the basis elements we obtain that
$$\beta_{21}=-\beta_{22},\quad \beta_{i1}=0,\ 3\leq i\leq n-1.$$
Now, considering $e_i=[e_{i-1},e_1]$ for $4\leq i\leq n,$ we have
\begin{equation*}
    \begin{split}
            D(e_i)&=D([e_{i-1},e_1])=[D(e_{i-1}),e_{1}]-[D(e_{1}),e_{i-1}]=0.
    \end{split}
\end{equation*}
The last relation completes the proof.
 \end{proof}

In following theorem, we give the description of the biderivations of $F_n^1$ and $F_n^2.$
\begin{thm}\label{F_n12}
    Any biderivations of $F_n^1$ and $F_n^2$ are described as follows:
    \begin{itemize}
        \item  for algebra $F_n^1:$
        \begin{equation*}
       \left\{ \left(
              \begin{array}{cccccccccc}
 \alpha_1   & 0       & 0      & \cdots & 0          & 0                  \\
 \alpha_2    & \alpha_1+\alpha_2   & 0 & \cdots & 0          & 0                \\
 \alpha_3    & \alpha_3   & 2\alpha_1+\alpha_2 & \cdots & 0          & 0                \\
 \cdots        &\cdots & \cdots        & \cdots & \cdots     & \cdots           \\
 \alpha_{n-1}& \alpha_{n-1}& \alpha_{n-2}& \cdots & (n-2)\alpha_1+\alpha_2 & 0            \\
 \alpha_n    &\alpha_{n+1} & \alpha_{n-1}&\cdots & \alpha_3 & (n-1)\alpha_1+\alpha_2        \\
 \end{array}
            \right),\ \left(
              \begin{array}{cccccccccc}
 \alpha_1    & 0             & \cdots & 0          & 0                  \\
 \beta_2    & 0  & \cdots & 0          & 0                \\
 \beta_3    & 0   & \cdots & 0          & 0                \\
 \cdots        & \cdots        & \cdots & \cdots     & \cdots           \\
 \beta_{n-1}& 0& \cdots & 0& 0            \\
 \beta_{n}    &\beta_{n+2} & \cdots & 0 & 0        \\
 \end{array}\right)\right\}
    \end{equation*}
        \item for algebra $F_n^2:$
    \begin{equation*}
       \left\{ \left(
              \begin{array}{cccccccccc}
 \alpha_1   & 0       & 0      & \cdots & 0          & 0                  \\
 \alpha_2    & \alpha_{n+1}   & 0 & \cdots & 0          & 0                \\
 \alpha_3    & 0   & 2\alpha_1 & \cdots & 0          & 0                \\
  \alpha_4    & 0   & \alpha_3 & \cdots & 0          & 0                \\
 \cdots        &\cdots & \cdots        & \cdots & \cdots     & \cdots           \\
 \alpha_{n-1}& 0& \alpha_{n-2}& \cdots & (n-2)\alpha_1 & 0            \\
 \alpha_n    &\alpha_{n+2} & \alpha_{n-1}&\cdots & \alpha_3 & (n-1)\alpha_1        \\
 \end{array}
            \right),\ \left(
              \begin{array}{cccccccccc}
 \alpha_1    & 0            &0 & \cdots & 0          & 0                  \\
 \beta_2    & \beta_{n+1}  &0& \cdots & 0          & 0                \\
 \beta_3    & 0   &0& \cdots & 0          & 0                \\
 \beta_4    & 0   &0& \cdots & 0          & 0                \\
 \cdots        & \cdots&\cdots        & \cdots & \cdots     & \cdots           \\
 \beta_{n-1}& 0& 0&\cdots & 0& 0            \\
 \beta_{n}    &\beta_{n+2} & 0& \cdots & 0 & 0        \\
 \end{array}\right)\right\}
    \end{equation*}
    \end{itemize}
\end{thm}
\begin{proof}
We prove the theorem for $F_n^1$, and one can do similar computations to obtain the description of biderivations of $F_n^2.$
Let $(d,D)$ be a biderivation of $NF_n,$ where $d\in Der(NF_n)$ and $D$ is an anti-derivation. We first apply the identity \eqref{1} for the elements $x=e_1, \ y=e_1$ and substitute the values from the derivation and anti-derivation. Then we have
    $$[e_1,\alpha_1e_1+\alpha_2e_2+\cdots+\alpha_ne_n]=[e_1,\beta_1e_1+\beta_2e_2+\cdots+\beta_ne_n]$$
Thus, we get $\beta_1=\alpha_1.$ Now, applying the same arguments for the elements $x=e_1, \ y=e_2,$ we have
    $$[e_1,(\alpha_1+\alpha_2)e_2+\alpha_3e_3\cdots+\alpha_{n-1}e_{n-1}+\alpha_{n+1}e_n]=[e_1,\beta_{n+1}e_1-\beta_{n+1}e_2+\beta_{n+2}e_n]$$
Hence, we get $\beta_{n+1}=0.$ We obtain the identities by checking the remaining products.
\end{proof}

Now, we describe the Leibniz algebra of ${\rm Bider}(F_n^1)$.
For this, let us fix the basis $\mathcal{B}$ of ${\rm Bider}(F_n^1)$. One can verify that ${\rm Bider}(F_n^1)$ is $(2n+1)$-dimensional. Thus, consider $\mathcal{B}=\{X_1,X_2,\dots,X_n,X_{n+1},Y_2,Y_3,\dots,Y_{n},Y_{n+1}\}$, where
\begin{align*}
    X_1&=\left(E_{11}+E_{22}+2E_{33}+\dots+(n-1)E_{nn},E_{11}\right)   & &\\
    X_2&=\left(E_{21}+E_{22}+\dots+E_{nn},0\right)    &  Y_2=&\left( 0,E_{21} \right)\\
     X_3&=\left(E_{31}+E_{32}+E_{43}+E_{54}+\dots+E_{nn-1},0\right)    &  Y_3=&\left( 0,E_{31} \right)\\
     X_4&=\left(E_{41}+E_{42}+E_{53}+\dots+E_{nn-2},0\right) &  Y_4=&\left( 0,E_{41} \right)\\
    &\dots & \dots\\
     X_{n-1}&=\left(E_{n-11}+E_{n-12}+E_{n3},0\right)   & Y_{n-1}=&\left( 0,E_{n-11} \right)\\
     X_{n}&=\left(E_{n1},0\right)   & Y_{n}=&\left( 0,E_{n1} \right),\\
      X_{n+1}&=\left(E_{n2},0\right)   & Y_{n+1}=&\left( 0,E_{n2} \right),
      \end{align*}
where $E_{ij}$ are matrix units.
Then it can be verified by \eqref{bider} that $[X_1,X_1]=[X_i,X_j]=0$ and $[X_1,Y_j]=[X_i,Y_j]=[Y_i,Y_j]=0,$ for $2\leq i,j\leq n$.
For other products we do the following multiplication table:
$$\begin{cases}
[X_1, X_2] = -Y_2,  & [X_{n+1}, X_1]= -(n-2)X_{n+1}, \\
[X_1, X_{n+1}]= (n-2)X_{n+1} - Y_{n+1}, & [Y_{n+1},X_1]=-(n-2)Y_{n+1},\\
[X_1, Y_2] = -Y_2, & [Y_2,X_{n+1}] =-Y_{2}, \\
[Y_k, X_2]=-Y_k, & \text{ for } 2 \leq k\leq n+1,\\
[Y_k, X_l] =-Y_{k+l-2}, & \text{ for } k-l+2 \leq n,\\
[Y_k,X_1] =-(k-2)Y_k, & \text{ for } 3 \leq k \leq n,\\
[X_k, X_1] = -(k-2)X_k, & \text{ for } 3 \leq k \leq n, \\
[X_1, Y_k] = (k-2)X_k - Y_k, & \text{  } 3 \leq k \leq n,
\end{cases}
$$
omitted products are zero.

Now, we describe the Leibniz algebra of ${\rm Bider}(F_n^2)$.
For this, let us fix the basis $\mathcal{B}$ of ${\rm Bider}(F_n^2)$. One can verify that ${\rm Bider}(F_n^2)$ is $(2n+3)$-dimensional. Thus, consider $\mathcal{B}=\{X_1,X_2,\dots,X_n,X_{n+1},X_{n+2},Y_2,Y_3,\dots,Y_{n+1},X_{n+2}\}$, where
\begin{align*}
    X_1&=\left(E_{11}+2E_{33}+\dots+(n-1)E_{nn},E_{11}\right)   & &\\
    X_2&=\left(E_{21},0\right)    &  Y_2=&\left( 0,E_{21} \right)\\
     X_3&=\left(E_{31}+E_{43}+E_{54}+\dots+E_{nn-1},0\right)    &  Y_3=&\left( 0,E_{31} \right)\\
     X_4&=\left(E_{41}+E_{53}+\dots+E_{nn-2},0\right) &  Y_4=&\left( 0,E_{41} \right)\\
    &\dots & \dots\\
     X_{n-1}&=\left(E_{n-11}+E_{n3},0\right)   & Y_{n-1}=&\left( 0,E_{n-11} \right)\\
     X_{n}&=\left(E_{n1},0\right)   & Y_{n}=&\left( 0,E_{n1} \right),\\
      X_{n+1}&=\left(E_{22},0\right)   & Y_{n+1}=&\left( 0,E_{22} \right),\\
      X_{n+2}&=\left(E_{n2},0\right)   & Y_{n+2}=&\left( 0,E_{n2} \right),
      \end{align*}
where $E_{ij}$ are matrix units.

Using equation \eqref{bider}, we obtain the following multiplication table with respect to the given basis:
$
\begin{cases}
    [X_1,X_2]=-Y_2,  & [X_2,X_1]=X_2, \\
    [X_1,X_n]=(n-2)X_n-Y_n & [X_n,X_1]=X_n  \\
    [Y_2,X_1]=Y_2. & [Y_{n+1},X_2]=Y_{2} \\
    [X_1,X_{n+2}]=(n-1)X_{n+2}, & [X_{n+2},X_1]=-(n-2)X_{n+2} \\
    [Y_{n+2},X_1]=-(n-1)Y_{n+2} &  [Y_{n+2},X_2]=Y_{n+2} \\
    [Y_{n+1},X_{n+2}]=Y_{n+2} & [Y_k,X_1]=-(k-2)Y_k, \text{ for  } 3\leq k\leq n,\\
    [X_1,X_k]=(k-2)X_k-Y_k\ \text{ for  } 3\leq k\leq n-1, & [X_k,X_1]=-(k-2)X_k,\ \text{ for  } 3\leq k\leq n-1, \\
    [Y_{k},X_l]=-Y_{k+l-2} & \text{ for } \ 3\leq l\leq n-1,\  3\leq k\leq n-1,\ k-l+2\leq n
\end{cases}$\\
omitted products are zero.

\section{Biderivations on solvable Leibniz algebras}
In this section we investigate derivations anti-derivations and biderivations on solvable Leibniz algebras with nilradicals $NF_n, F_n^1$ and $F_n^2.$

\subsection{Biderivations on solvable Leibniz algebras with nilradicals $NF_n$}

Now, we present the classification result for arbitrary dimensional solvable
Leibniz algebras with null-filiform nilradical.
\begin{thm}[\cite{CLOK12}]
Let $R$ be a solvable Leibniz algebra whose nilradical is $NF_n.$ Then there exists a basis $\{e_1,e_2,...,e_n,h\}$ of the algebra $R$ such that the multiplication table of $R$ with respect to this basis has the following from:
  \begin{align*}
  [e_i,e_1]&=e_{i+1}, \quad 1\leq i\leq n-1,\\
  [h,e_1]&=-e_{1},\\
  [e_i,h]&=ie_i, \quad 1\leq i\leq n
\end{align*}
\end{thm}

Now we give the description of a biderivations on the solvable Leibniz algebra $R.$ This is done by first analyzing its derivations and then its anti-derivations.

\begin{prop}\cite{der}
Let $R$ be a solvable Leibniz algebra whose nilradical is $NF_n.$ Any derivation of the algebra $R$ has the following form:
\begin{equation*}
 \left(\begin{array}{cccccccccc}
 0   & 0       & 0      & \cdots & 0          & 0                  \\
 -\alpha_1    & \alpha_2   & 0 & \cdots & 0          & 0                \\
 0    & \alpha_1  & 2\alpha_2 & \cdots & 0          & 0                \\
  0    & 0 & \alpha_1 & \cdots & 0          & 0                \\
 \cdots        &\cdots & \cdots        & \cdots & \cdots     & \cdots           \\
 0& 0& 0& \cdots & (n-1)\alpha_2 & 0            \\
 0    &0& 0&0 & \alpha_1 & n\alpha_2        \\
 \end{array}
            \right)
    \end{equation*}
\end{prop}

\begin{prop}
Any anti-derivation of $R$ has the following form:
\begin{equation*}
\left(\begin{array}{cccccccccc}
 0   & 0       & 0      & \cdots & 0          & 0                  \\
 \beta_2    & \beta_1   & 0 & \cdots & 0          & 0                \\
 2\beta_3    & \beta_2  & 0 & \cdots & 0          & 0                \\
  3\beta_4    & \beta_3  & 0 & \cdots & 0          & 0                \\
 \cdots        &\cdots & \cdots        & \cdots & \cdots     & \cdots           \\
 (n-1)\beta_{n}& \beta_{n-1}& 0& \cdots & 0 & 0            \\
 \beta_{n+1}    &\beta_n & 0&0 & 0 & 0        \\
 \end{array}
            \right)
    \end{equation*}
\end{prop}

\begin{proof}
     It is easy to see that  from the multiplication table of the algebra $R$ we conclude that $e_1$ and $h$ are the generator basis elements of the algebra. We use these generators to calculate the anti-derivation. Let $D$ be an anti-derivations of $R.$ Set
$$D(h)=\sum\limits_{i=2}^{n}\beta_{1i}e_i+\beta_{11}h, \quad
D(e_1)=\sum\limits_{i=2}^{n}\beta_{2i}e_i+\beta_{2 1}h.$$
Since $D([x,x])=0,$ we can easily obtain that $D(e_2)=0.$ Now, by $e_i=[e_{i-1},e_1],$ we have
\begin{equation*}
  \begin{split}
      D(e_i)&=[D(e_{i-1}),e_1]-[D(e_1),e_{i-1}]=0.\quad 3\leq i\leq n.
  \end{split}
\end{equation*}
Moreover, considering $[h,e_1]=-e_1,$ we have the following chain of equalities:
\begin{equation*}
    \begin{split}
       -\sum\limits_{i=2}^{n}\beta_{2i}e_i-\beta_{2 1}h= -D(e_1)&=[D(h),e_1]-[D(e_1),h]=\sum\limits_{i=2}^n\beta_{1i}[e_i,e_1]+\beta_{1 1}[h,e_1]-\sum\limits_{i=2}^{n}\beta_{2i}[e_i,h]-\beta_{2 1}[h,h]\\
        &=\sum\limits_{i=2}^n\beta_{1  i-1}e_i-\beta_{11}e_1-\sum\limits_{i=1}^{n}i\beta_{2i}e_i.
    \end{split}
\end{equation*}
Comparing coefficients of the basis elements we obtain that
$$\beta_{21}=\beta_{11}=0,\quad \beta_{1  i-1}=(i-1)\beta_{2i},\quad 2\leq i\leq n. $$
  The remaining products reduce to an identity.
\end{proof}

Now, we are ready to give the description of the biderivations of $R.$
\begin{thm}
Any biderivations of $R$ is described as follows:
    \begin{equation}
       \left\{  \left(
              \begin{array}{cccccccccc}
 0   & 0       & 0      & \cdots & 0          & 0                  \\
 -\alpha_1    & \alpha_2   & 0 & \cdots & 0          & 0                \\
 0    & \alpha_1  & 2\alpha_2 & \cdots & 0          & 0                \\
  0    & 0 & \alpha_1 & \cdots & 0          & 0                \\
 \cdots        &\cdots & \cdots        & \cdots & \cdots     & \cdots           \\
 0& 0& 0& \cdots & (n-1)\alpha_2 & 0            \\
 0    &0& 0&0 & \alpha_1 & n\alpha_2        \\
 \end{array}
            \right),\ \left(
              \begin{array}{cccccccccc}
 0   & 0       & 0      & \cdots & 0          & 0                  \\
 -\alpha_1    & \alpha_2   & 0 & \cdots & 0          & 0                \\
 2\beta_3    & -\alpha_1  & 0 & \cdots & 0          & 0                \\
  3\beta_4    & \beta_3  & 0 & \cdots & 0          & 0                \\
 \cdots        &\cdots & \cdots        & \cdots & \cdots     & \cdots           \\
 (n-1)\beta_{n}& \beta_{n-1}& 0& \cdots & 0 & 0            \\
 \beta_{n+1}    &\beta_n & 0&0 & 0 & 0        \\
 \end{array}
            \right)\right\}
    \end{equation}
\end{thm}
\begin{proof}
The proof of the theorem is done similarly as the proofs of the previous theorems by applying the equality \eqref{1} to all possible pairs of elements $x,y\in R.$
\end{proof}

\begin{rem}
  It should be noted that any biderivation defined on $R$ is an inner biderivation.
\end{rem}
Now, we describe the Leibniz algebra ${\rm Bider}(R)$.
For this, let us fix the basis $\mathcal{B}$ of ${\rm Bider}(R)$. One can verify that ${\rm Bider}(R)$ is $(n+2)$-dimensional. Thus, consider $\mathcal{B}=\{X_1,X_2,X_3,Y_4,Y_5,\dots, Y_{n+2}\}$, where
\begin{center}
$H:= \bigg( E_{22}+2E_{33}+3E_{44}+ \ldots +n E_{n+1 \ n+1} , E_{22} \bigg);$

$X_1:= \bigg( -E_{21}+E_{32}+E_{43}+ \ldots +E_{n+1 \ n} , -E_{21}-E_{32} \bigg);$ 

$X_{i}:=\bigg(0 , i E_{i+1\, 1}+E_{i+2 \, 2} \bigg), \ i \in \lbrace 2, 3, \ldots n-1 \rbrace ;$ \quad $X_{n}:= \bigg(0 ,E_{n+1 \ 1}\bigg)$
\end{center}
where $E_{ij}$ are matrix units.

Based on equation \eqref{bider}, the multiplication table corresponding to the above basis takes the following form:
$$
\begin{cases}
\big[H,X_1 \big] = X_1, &   \\
\big[X_{i} , X_1 \big] = -X_{i+1},& \big[X_i , H \big] = i X_{i},\\
\end{cases}
$$
where $i \in \lbrace 1,3 \ldots n \rbrace.$ Omitted products are zero.

\subsection{Biderivations on solvable Leibniz algebras with naturally graded filiform non-Lie Leibniz algebra}
Solvable Leibniz algebras whose nilradical is the naturally graded filiform Leibniz algebra $F_n^1$ are classified in \cite{ladra}. It is proved that the dimension of a solvable Leibniz algebra whose nilradical is isomorphic to an $n$-dimensional naturally graded filiform Leibniz algebra is not greater than $n+2$.

\begin{thm} Let $R(F_n^1)$ be $(n+2)$-dimensional solvable Leibniz algebra with nilradical $F_n^1$. Then it is siomorphic to the Leibniz algebra with a basis
$\{e_1,e_2,\dots,e_n,h_1,h_2\}$ and the multiplication table:
\[\begin{array}{llll}
[e_i,e_1]=e_{i+1}, &  2\leq i \leq {n-1},  & [e_1,h_2]=e_1,      &                  \\[1mm]
[e_i,h_1]=e_i,       &  2\leq i\leq n,       & [e_i,h_2]=(i-1)e_i, &  2\leq i \leq n, \\[1mm]
                   &                       & [h_2,e_1]=-e_1.     &
\end{array}\]
\end{thm}

Next, we present a description of the biderivations of the algebra  $R(F_n^1)$. This is done by first analyzing its derivations and then its anti-derivations.
We first present the classification of derivations of the algebra $R(F_n^1)$ as given in \cite{KhO}.
\begin{prop}\cite{KhO}
Let $R(F_n^1)$ be $(n+2)$-dimensional solvable Leibniz algebra with nilradical $F_n^1$. Any derivation of the algebra $R(F_n^1)$ has the following form:

\begin{equation*}
     \left(\begin{array}{cccccccccc}
 0  &0 & 0       & 0      & \cdots & 0          & 0                  \\
  0 &0 & 0       & 0      & \cdots & 0          & 0                  \\
 0 &-\alpha_2   & \alpha_3   & 0 & \cdots & 0          & 0                \\
 0  &0  & 0  & \alpha_1+\alpha_3 & \cdots & 0          & 0                \\
  0  &0  & 0 & \alpha_2 & \cdots & 0          & 0                \\
 \cdots& \cdots        &\cdots & \cdots        & \cdots & \cdots     & \cdots           \\
 0& 0&0& 0& \cdots & \alpha_1+(n-2)\alpha_3 & 0            \\
 0 &0   &0& 0&0 & \alpha_2 & \alpha_1+(n-1)\alpha_3        \\
 \end{array}
            \right)
    \end{equation*}
\end{prop}

Now, we compute anti-derivation of $R(F_n^1)$

\begin{prop}
Let $R(F_n^1)$ be $(n+2)$-dimensional solvable Leibniz algebra with nilradical $F_n^1$. Any anti-derivation of the algebra $R(F_n^1)$ has the following form:
\begin{equation*}
   \left(\begin{array}{cccccccccc}
 \beta_{n+2}  &\beta_{n+3} & 0       & 0      & \cdots & 0          & 0                  \\
  0 &0 & 0       & 0      & \cdots & 0          & 0                  \\
 0 &\beta_2   & \beta_1   & 0 & \cdots & 0          & 0                \\
 \beta_3  &\beta_3  & 0  & 0 & \cdots & 0          & 0                \\
  \beta_4  &2\beta_4  & \beta_3 & 0 & \cdots & 0          & 0                \\
 \cdots& \cdots        &\cdots & \cdots        & \cdots & \cdots     & \cdots           \\
 \beta_n& (n-2)\beta_n&\beta_{n-1}& 0& \cdots & 0 & 0            \\
 \beta_{n+1} & (n-1)\beta_{n+1}   &\beta_n& 0&0 & 0 & 0        \\
 \end{array}
 \right)
\end{equation*}
\end{prop}

\begin{proof}
Let $D$ be an anti-derivations of $R(F_n^1).$ We first set
$$D(e_1)=\sum\limits_{i=1}^{n}\alpha_ie_i+\gamma_1 h_1+\beta_1 h_2,\quad
D(e_2)=\sum\limits_{i=1}^{n}\nu_ie_i+\gamma_2h_1+\beta_2h_2,$$
$$
D(h_1)=\sum\limits_{i=1}^{n}\mu_ie_i+\gamma_4h_1+\beta_4h_2,\quad
D(h_2)=\sum\limits_{i=1}^{n}\lambda_ie_i+\gamma_3h_1+\beta_3h_2.$$
Now we consider $e_2=[e_2,h_1]$ and apply the identity of anti-derivation.
\begin{equation*}
    \begin{split}
      \sum\limits_{i=1}^{n}\nu_ie_i+\gamma_2h_1+\beta_2h_2=  D(e_2)&=[D(e_2),h_1]-[D(h_1),e_2]=\sum\limits_{i=1}^{n}\nu_i[e_i,h_1]=\sum\limits_{i=2}^{n}\nu_ie_i.
    \end{split}
\end{equation*}
Then we have $\nu_1=\gamma_2=0.$ Now considering $e_3=[e_2,e_1]$ and $2e_3=[e_3,h_2],$ we have the following:
\begin{equation*}
  \begin{split}
      D(e_3)&=[D(e_2),e_1]-[D(e_1),e_2]=\sum\limits_{i=2}^{n}\nu_i[e_i,e_1]+\beta_2[h_2,e_1]=\sum\limits_{i=3}^{n}\nu_{i-1}e_i-\beta_2e_1,\\
     2 D(e_3)&=[D(e_{3}),h_2]-[D(h_2),e_{3}]=\sum\limits_{i=3}^{n}\nu_{i-1}[e_i,h_2]-\beta_2[e_1,h_2]=\sum\limits_{i=3}^{n}(i-1)\nu_{i-1}e_i-\beta_2e_1,
  \end{split}
\end{equation*}
which imply that $\beta_2=\nu_i=0,\ 3\leq i\leq n-1.$
Similarly, by $e_1=[e_1,h_2],$ we get
\begin{equation*}
    \begin{split}
      \sum\limits_{i=1}^{n}\alpha_ie_i+\gamma_1h_1+\beta_1h_2=  D(e_1)&=[D(e_1),h_2]-[D(h_2),e_1]=\sum\limits_{i=1}^{n}\alpha_i[e_i,h_2]-\sum\limits_{i=1}^{n}\lambda_i[e_i,e_1]-\beta_3[h_2,e_1]\\
      &=\alpha_1e_1+\sum\limits_{i=2}^{n}(i-1)\alpha_ie_i-\sum\limits_{i=3}^{n}\lambda_{i-1}e_i+\beta_3e_1
    \end{split}
\end{equation*}
This follows that $\beta_1=\beta_3=\gamma_1=0,\ \lambda_{i-1}=(i-2)\alpha_i,\ 3\leq i\leq n.$
Moreover, we have $[e_1,h_1]=0.$ Then applying the definition of anti-derivation for the pair $e_1$ and $h_1$ we have
\begin{equation*}
    \begin{split}
      0&=D([e_1,h_1])=[D(e_1),h_1]-[D(h_1),e_1]=\sum\limits_{i=1}^{n}\alpha_i[e_i,h_1]-\sum\limits_{i=1}^{n}\mu_i[e_i,e_1]-\beta_4[h_2,e_1]\\
      &=\sum\limits_{i=2}^{n}\alpha_ie_i-\sum\limits_{i=3}^{n}\mu_{i-1}e_i+\beta_4e_1.
    \end{split}
\end{equation*}
which leads to the restrictions $\alpha_2=\beta_4=0,\ \mu_{i-1}=\alpha_i,\ 3\leq i\leq n.$

Similarly, taking $[h_2,e_2]=0$ into account, we have the following:
\begin{equation*}
    \begin{split}
      0&=D([h_2,e_2])=[D(h_2),e_2]-[D(e_2),h_2]=-\nu_2e_2-(n-1)\nu_ne_n.
    \end{split}
\end{equation*}
Then we obtain that $\nu_2=\nu_n=0.$ Similarly, we have the following chains of equations:
\begin{equation*}
  \begin{split}
      D(e_3)&=[D(e_2),e_1]-[D(e_1),e_2]=0,\\
      D(e_i)&=[D(e_{i-1}),e_1]-[D(e_1),e_{i-1}]=0.\quad 4\leq i\leq n.
  \end{split}
\end{equation*}

\begin{equation*}
  \begin{split}
      0&=[D(h_2),h_1]-[D(h_1),h_2]\\
      &=\lambda_1[e_1,h_1]+\sum\limits_{i=2}^{n-1}(i-1)\alpha_{i+1}[e_i,h_1]+\lambda_n [e_n,h_1]-\mu_1[e_1,h_2]-\sum\limits_{i=2}^{n-1}\alpha_{i+1}[e_i,h_2]-\mu_n [e_n,h_2]\\
      &=-\mu_{1}e_1+\lambda_ne_n-(n-1)\mu_ne_n
  \end{split}
\end{equation*}
These imply that $\mu_1=0,\ \lambda_n=(n-1)\mu_n.$
  The remaining products reduce to an identity. Summarizing all the relations obtained above, we get the required result.
\end{proof}

Based on the previously computed derivations and anti-derivations, we describe the biderivations in what follows.

\begin{thm}
    Any biderivations of $R(F_n^1)$  is described as follows:
{\small
\begin{equation*}
     {  \left\{ \left(
              \begin{array}{cccccccccc}
 0  &0 & 0       & 0      & \cdots & 0          & 0                  \\
  0 &0 & 0       & 0      & \cdots & 0          & 0                  \\
 0  &-\alpha_2    & \alpha_3   & 0 & \cdots & 0          & 0                \\
 0  &0  & 0  & \alpha_1+\alpha_3 & \cdots & 0          & 0                \\
  0  &0  & 0 & \alpha_2 & \cdots & 0          & 0                \\
 \cdots& \cdots        &\cdots & \cdots        & \cdots & \cdots     & \cdots           \\
 0& 0&0& 0& \cdots & \alpha_1+(n-2)\alpha_3 & 0            \\
 0 &0   &0& 0&0 & \alpha_2 & \alpha_1+(n-1)\alpha_3        \\
 \end{array}
\right),\ \left(
\begin{array}{cccccccccc}
  0 &0 & 0       & 0      & \cdots & 0                            \\
  0 &0 & 0       & 0      & \cdots & 0                            \\
  0 &-\alpha_2   & \alpha_3   & 0 & \cdots & 0                          \\
 \beta_3  &\beta_3  & 0  & 0 & \cdots           & 0                \\
 \beta_4  &2\beta_4  & \beta_3 & 0 & \cdots           & 0                \\
 \cdots& \cdots        &\cdots & \cdots         & \cdots     & \cdots           \\
 \beta_n& (n-2)\beta_n&\beta_{n-1}& 0& \cdots & 0             \\
 \beta_{n+1} &(n-1)\beta_{n+1}   &\beta_n& 0&\cdots & 0         \\
 \end{array}
            \right)\right\}}
\end{equation*}
}\end{thm}
\begin{proof}
    Let $(d,D)$ be a biderivation of $R(F_n^1).$ Then by
    $[e_2,d(e_1)]=[e_2,D(e_1)],$
    we have
    $$[e_1,\alpha_3 e_1]=[e_1,\beta_1e_1+\beta_3e_3+\cdots+\beta_ne_n],$$
    which follows that $\beta_1=\alpha_3.$
   Similarly, by
    $[e_2,d(h_2)]=[e_2,D(h_2)]$ and substituting the values from the derivation and anti-derivation, we can write the following equality:
    $$[e_2,-\alpha_2e_1]=[e_2,\beta_{n+3} h_1+\beta_2 e_1+\beta_3e_2+\cdots+(n-2)\beta_ne_{n-1}+(n-1)\beta_{n+1} e_n],$$
    which leads to the restrictions $\alpha_2=-\beta_2,\ \beta_{n+3}=0.$
Applying the same arguments for the elements $x=e_2, \ y=h_1,$ we have and substituting the values from the derivation and anti-derivation, we have
    $$0=[e_2,\beta_{n+2} h_1+\beta_3e_2+\cdots+\beta_ne_{n-1}+\beta_{n+1} e_n].$$
    This follows that $\beta_{n+2}=0.$ Note that we obtain the identities by checking the remaining products.
   \end{proof}
Now, we describe the Leibniz algebra ${\rm Bider}(R(F_n^1))$.
For this, let us fix the basis $\mathcal{B}$ of ${\rm Bider}(R(F_n^1))$. One can verify that ${\rm Bider}(R(F_n^1))$ is $(n+2)$-dimensional. Thus, consider $\mathcal{B}=\{X_1,X_2,X_3,Y_4,Y_5,\dots, Y_{n+2}\}$, where
\begin{center}
$
H_1:=-\big( E_{44}+E_{55}+ \ldots+E_{n+2 \ n+2} , 0 \big),
$

$
H_2:=\big( E_{33}+E_{44}+2E_{55}+ \ldots+(n-1)E_{n+2 \ n+2} , E_{33}\big),
$

$
X_1:=-\big( -E_{32}+E_{54}+E_{65}+ \ldots+E_{n+2 \ n+1} , -E_{32}\big),
$

$
X_{i}:=\big( 0 , E_{i+2 \ 1}+(i-1)E_{i+2 \ 2}+E_{i+3 \ 3}\big), i \in \lbrace 2,3, \ldots n-1 \rbrace
$

$
X_{n}:=\big( 0 , E_{n+2 \ 1}+(n-1)E_{n+2 \ 2}\big),
$
\end{center}
where $E_{ij}$ are matrix units.

By applying equation \eqref{bider}, we derive the following multiplication table relative to the basis given above:
$$
\begin{cases}
\big[ X_1, H_2 \big]=X_1, & \big[ H_2, X_1 \big]=-X_1,  \\
\big[ X_i, H_1 \big]=X_{i} & (i \in \lbrace 2,3, \ldots n \rbrace), \\
\big[ X_i, H_2 \big]=(i-1)X_{i} & (i \in \lbrace 2,3, \ldots n \rbrace), \\
\big[ X_i, X_1 \big]=X_{i+1} & (i \in \lbrace 2,3, \ldots n-1 \rbrace).
\end{cases}
$$
omitted products are zero.

Solvable Leibniz algebras whose nilradical is the naturally graded filiform Leibniz algebra $F_n^2$ are classified in \cite{CLOK13}. It is proved that the dimension of a solvable Leibniz algebra whose nilradical is isomorphic to an $n$-dimensional naturally graded filiform Leibniz algebra is not greater than $n+2$.

\begin{thm}\label{thm5.8} An arbitrary $(n+2)$-dimensional solvable Leibniz algebra with nilradical $F_n^2$
 is isomorphic to one of the following non-isomorphic algebras:
\[ \mathcal{L}_1: \left\{\begin{aligned}
{} [e_1,e_1] & =e_3, & [e_i,e_1] & =e_{i+1}, && 3\leq i\leq n-1,\\
 [e_1,h_2]& =e_1, & [h_2,e_1] & =-e_1,   && \\
 [e_2,h_1]& =-[h_1,e_2]=e_2, &[e_i,h_2]&=(i-1)e_i, && 3\leq i\leq n,
\end{aligned}\right.\]
\[ \mathcal{L}_2: \left\{\begin{aligned}
{} [e_1,e_1] & =e_3, & [e_i,e_1] & =e_{i+1}, && 3\leq i\leq n-1,\\
 [e_1,h_2]& =e_1, & [h_2,e_1] & =-e_1,   && \\
 [e_2,h_1]& =e_2, &[e_i,h_2]&=(i-1)e_i, && 3\leq i\leq n\,.
\end{aligned}\right.\]
\end{thm}

We first describe the biderivations of $\mathcal{L}_1$ and $\mathcal{L}_2$ by considering their derivations and then their anti-derivations. In what follows, we provide the classification of derivations of the algebras $\mathcal{L}_1$ and $\mathcal{L}_2$ , as described in \cite{KhS}.
\begin{prop}
Let $(n+2)$-dimensional solvable Leibniz algebra with nilradical $F_n^2$. Any derivation of the algebras $\mathcal{L}_1$ and $\mathcal{L}_2$ has the following form:
\begin{itemize}
    \item For algebra $\mathcal{L}_1:$
    \begin{equation*}
       \left(
              \begin{array}{cccccccccc}
 0  &0 & 0       & 0      &0& \cdots & 0          & 0                  \\
  0 &0 & 0       & 0      &0& \cdots & 0          & 0                  \\
 0 &-\alpha_1   & \alpha_4   & 0 & 0&\cdots & 0          & 0                \\
-\alpha_2 &0 & 0 & \alpha_3 & 0&\cdots & 0          & 0                \\
   0 &0  & \alpha_1 & 0&2\alpha_4 &\cdots & 0          & 0                \\
   0 &0  &0& 0&\alpha_1&\cdots&0 & 0        \\
 \cdots& \cdots      &\cdots  &\cdots & \cdots        & \cdots & \cdots     & \cdots           \\
 0&0& 0&0& 0&\cdots & (n-2)\alpha_4 &0           \\
 0 &0  &0&0&0 &\cdots & \alpha_1 & (n-1)\alpha_4       \\
 \end{array}
            \right)
\end{equation*}
    \item For algebra $\mathcal{L}_2:$
    \begin{equation*}
       \left(
              \begin{array}{cccccccccc}
 0  &0 & 0       & 0      &0& \cdots & 0          & 0                  \\
  0 &0 & 0       & 0      &0& \cdots & 0          & 0                  \\
 0 &-\alpha_1   & \alpha_3   & 0 & 0&\cdots & 0          & 0                \\
0 &0 & 0 & \alpha_2 & 0&\cdots & 0          & 0                \\
   0 &0  & \alpha_1 & 0&2\alpha_3 &\cdots & 0          & 0                \\
   0 &0  &0& 0&\alpha_1&\cdots&0 & 0        \\
 \cdots& \cdots      &\cdots  &\cdots & \cdots        & \cdots & \cdots     & \cdots           \\
 0&0& 0&0& 0&\cdots & (n-2)\alpha_3 &0           \\
 0 &0  &0&0&0 &\cdots & \alpha_1 & (n-1)\alpha_3       \\
 \end{array}
            \right)
\end{equation*}

\end{itemize}
\end{prop}

Below, we present the description of the spaces of anti-derivations of $\mathcal{L}_1$ and $\mathcal{L}_2$ .

\begin{prop}
Let $(n+2)$-dimensional solvable Leibniz algebra with nilradical $F_n^2$. Any anti-derivation of the algebras $\mathcal{L}_1$ and $\mathcal{L}_2$ has the following form:
\begin{itemize}
    \item For algebra $\mathcal{L}_1:$
    \begin{equation*}
       \left(
              \begin{array}{cccccccccc}
 0  &0 & 0       & 0      & \cdots & 0          & 0                  \\
  0 &0 & 0       & 0      & \cdots & 0          & 0                  \\
 0 &\beta_3   & \beta_1   & 0 & \cdots & 0          & 0                \\
   \beta_{n+2}  &0 & 0 & \beta_2 & \cdots & 0          & 0                \\
   0 &2\beta_4  & \beta_3 & 0 & \cdots & 0          & 0                \\
 \cdots& \cdots        &\cdots & \cdots        & \cdots & \cdots     & \cdots           \\
 0& (n-2)\beta_n&\beta_{n-1}& 0& \cdots & 0 & 0            \\
 0 &\beta_{n+1}   &\beta_n& 0&0 & 0 & 0        \\
 \end{array}
            \right).
\end{equation*}
    \item For algebra $\mathcal{L}_2:$
    \begin{equation*}
       \left(
              \begin{array}{cccccccccc}
 \beta_{n+2}  &\beta_2 & 0       & 0      & \cdots & 0          & 0                  \\
  0 &0 & 0       & 0      & \cdots & 0          & 0                  \\
 0 &\beta_3   & \beta_1   & 0 & \cdots & 0          & 0                \\
   \beta_{n+3}  &0 & 0 & 0 & \cdots & 0          & 0                \\
   0 &2\beta_4  & \beta_3 & 0 & \cdots & 0          & 0                \\
 \cdots& \cdots        &\cdots & \cdots        & \cdots & \cdots     & \cdots           \\
 0& (n-2)\beta_n&\beta_{n-1}& 0& \cdots & 0 & 0            \\
 0 &\beta_{n+1}   &\beta_n& 0&0 & 0 & 0        \\
 \end{array}
            \right).
\end{equation*}

\end{itemize}

\end{prop}

Using the derivations and anti-derivations obtained above, we now provide the description of the biderivations below.

\begin{thm}
    Any biderivations of $\mathcal{L}_1$ and $\mathcal{L}_2$   $(n+2)$-dimensional solvable Leibniz algebra with nilradical $F_n^2$ can be described as follows:
{\small
\begin{equation*}
     {\left\{ \left(
              \begin{array}{cccccccccc}
 0  &0 & 0       & 0      &0& \cdots & 0          & 0                  \\
  0 &0 & 0       & 0      &0& \cdots & 0          & 0                  \\
 0 &-\alpha_1   & \alpha_4   & 0 & 0&\cdots & 0          & 0                \\
-\alpha_2 &0 & 0 & \alpha_3 & 0&\cdots & 0          & 0                \\
   0 &0  & \alpha_1 & 0&2\alpha_4 &\cdots & 0          & 0                \\
   0 &0  &0& 0&\alpha_1&\cdots&0 & 0        \\
 \cdots& \cdots      &\cdots  &\cdots & \cdots        & \cdots & \cdots     & \cdots           \\
 0&0& 0&0& 0&\cdots & (n-2)\alpha_4 &0           \\
 0 &0  &0&0&0 &\cdots & \alpha_1 & (n-1)\alpha_4       \\
 \end{array}
            \right),\ \left(
              \begin{array}{cccccccccc}
 0  &0 & 0       & 0      & \cdots & 0          & 0                  \\
  0 &0 & 0       & 0      & \cdots & 0          & 0                  \\
 0 &-\alpha_1   & \alpha_4  & 0 & \cdots & 0          & 0                \\
   -\alpha_2  &0 & 0 & \alpha_3 & \cdots & 0          & 0                \\
   0 &2\beta_4  & -\beta_1 & 0 & \cdots & 0          & 0                \\
 \cdots& \cdots        &\cdots & \cdots        & \cdots & \cdots     & \cdots           \\
 0& (n-2)\beta_n&\beta_{n-1}& 0& \cdots & 0 & 0            \\
 0 &\beta_{n+1}   &\beta_n& 0&\cdots & 0 & 0        \\
 \end{array}
            \right)\right\}}
\end{equation*}

{\small \begin{equation*}
     { \left\{  \left(
              \begin{array}{cccccccccc}
 0  &0 & 0       & 0      &0& \cdots & 0          & 0                  \\
  0 &0 & 0       & 0      &0& \cdots & 0          & 0                  \\
 0 &-\alpha_1   & \alpha_3   & 0 & 0&\cdots & 0          & 0                \\
0 &0 & 0 & \alpha_2 & 0&\cdots & 0          & 0                \\
   0 &0  & \alpha_1 & 0&2\alpha_3 &\cdots & 0          & 0                \\
   0 &0  &0& 0&\alpha_1&\cdots&0 & 0        \\
 \cdots& \cdots      &\cdots  &\cdots & \cdots        & \cdots & \cdots     & \cdots           \\
 0&0& 0&0& 0&\cdots & (n-2)\alpha_3 &0           \\
 0 &0  &0&0&0 &\cdots & \alpha_1 & (n-1)\alpha_3       \\
 \end{array}
            \right),\  \left(
              \begin{array}{cccccccccc}
0  &0& 0       & 0      & \cdots & 0          & 0                  \\
0 &0 & 0       & 0      & \cdots & 0          & 0                  \\
0 &-\alpha_1   & \alpha_3  & 0 & \cdots & 0          & 0                \\
\beta_{n+2}  &0 & 0 & 0 & \cdots & 0          & 0                \\
0 &2\beta_4  & -\alpha_1 & 0 & \cdots & 0          & 0                \\
\cdots& \cdots        &\cdots & \cdots        & \cdots & \cdots     & \cdots           \\
0& (n-2)\beta_n&\beta_{n-1}& 0& \cdots & 0 & 0            \\
0 &\beta_{n+1}   &\beta_n& 0&\cdots & 0 & 0        \\
 \end{array}
\right)\right\}}
    \end{equation*}}}
\end{thm}
 \begin{proof}
 We skip the proof of the theorem as it is similar to the proof of the previous theorems.
 \end{proof}
Now, we describe the Leibniz algebra ${\rm Bider}(\mathcal{L}_1)$.
For this, let us fix the basis $\mathcal{B}$ of ${\rm Bider}(\mathcal{L}_1)$. One can verify that ${\rm Bider}(\mathcal{L}_1)$ is $(n+2)$-dimensional. Thus, consider $\mathcal{B}=\{X_1,X_2,X_3,X_4,Y_5,\dots, Y_{n+2}\}$, where
\begin{center}
{$
H_1:=\big( E_{44} , E_{44} \big), \quad
H_2:=\big( E_{33}+2E_{55}+3E_{66}+ \ldots (n-1)E_{n+2 \ n+2},E_{33}\big),
$

$
X_1:=\big( -E_{32}+E_{53}+E_{65}+ \ldots+E_{n+2 \ n+1} , -E_{32}-E_{53}\big), \quad X_2:=\big( -E_{41} , -E_{41}\big),
$

$
X_{i}:=\big( 0 , (i-1)E_{i+2 \ 2}+E_{i+3 \ 3}\big), i \in \lbrace 3,4  \ldots n-1 \rbrace \ , \quad X_{n}:=\big( 0 , E_{n+2 \ 2}\big)
$}
\end{center}
where $E_{ij}$ are matrix units.
Then we have the following multiplication table with respect to the basis introduced above:

\textcolor{black}{$\begin{cases}
\big[ X_1, X_1 \big]=X_3 \,,& \big[ X_1, X_4 \big]=-\big[ X_4, X_1 \big]=-X_1, \\
\big[ X_2, H_1 \big]=-\big[ H_1, X_2 \big]=-X_2,  &  \big[ X_i, X_1 \big]=-X_{i+1}(i \in \lbrace 3,4, \ldots n-1 \rbrace),\\
\big[ X_i, H_2 \big]=-(i-1)X_{i}(i \in \lbrace 3,4, \ldots n-1 \rbrace),  & \big[ X_{n}, H_2 \big]=-(n-1)X_{n}.
\end{cases}$}

omitted products are zero.

Now, we describe the Leibniz algebra ${\rm Bider}(\mathcal{L}_2)$.
For this, let us fix the basis $\mathcal{B}$ of ${\rm Bider}(\mathcal{L}_2)$. One can verify that ${\rm Bider}(\mathcal{L}_2)$ is $(n+2)$-dimensional. Thus, consider $\mathcal{B}=\{X_1,X_2,X_3,X_4,Y_5,\dots, Y_{n+2}\}$, where

\begin{center}
{$
 H_1:=\big( E_{44} , 0 \big), \quad
H_2:=\big( E_{33}+2E_{55}+3E_{66}+ \ldots (n-1)E_{n+2 \ n+2},E_{33}\big),
$

$
X_1:=\big( -E_{32}+E_{53}+E_{65}+ \ldots+E_{n+2 \ n+1} , -E_{32}-E_{53}\big), \quad X_2:=\big( 0 , E_{41}\big),
$

$
X_{i}:=\big( 0 , (i-1)E_{i+2 \ 2}+E_{i+3 \ 3}\big), i \in \lbrace 3,4  \ldots n-1 \rbrace \ , \quad X_{n}:=\big( 0 , E_{n+2 \ 2}\big)
$}
\end{center}

where $E_{ij}$ are matrix units.

Using equation \eqref{bider}, we obtain the following multiplication table with respect to the given basis:

\textcolor{black}{$\begin{cases}
\big[ X_1, X_1 \big]=X_3 \,,& \big[ X_1, X_4 \big]=-\big[ X_4, X_1 \big]=-X_1, \\
\big[ X_2, H_1 \big]=-X_2,  &  \big[ X_i, X_1 \big]=-X_{i+1}(i \in \lbrace 3,4, \ldots n-1 \rbrace),\\
\big[ X_i, H_2 \big]=-(i-1)X_{i}(i \in \lbrace 3,4, \ldots n-1 \rbrace),  & \big[ X_{n}, H_2 \big]=-(n-1)X_{n}.
\end{cases}$}

omitted products are zero.

\begin{rem}
It should be mentioned that the Leibniz algebras of biderivations of all algebras considered in this section are isomorphic to the corresponding Leibniz algebras under consideration.
\end{rem}

\section{Acknowledgements}
We thank professor A.Kh.Khudoyberdiyev for the helpful comments and suggestions that contributed
to improving this paper.

\end{document}